\documentclass[11pt]{article}

\usepackage{amsmath,amssymb,psfig,epsfig,latexsym,graphicx,here}

\newcommand{\bsq}{{\vrule height .9ex width .8ex depth -.1ex }}
\newcommand{\eeq}{\end{equation}}
\newcommand{\beql}[1]{\begin{equation}\label{#1}}
\newcommand{\beq}{\begin{equation}}
\newcommand{\eqn}[1]{(\ref{#1})}
\newcommand{\ZZ}{{\mathbb Z}}
\newcommand{\RR}{{\mathbb R}}

\newcommand{\sQ}{{\mathcal Q}}
\newcommand{\sO}{{\mathcal O}}
\newcommand{\sP}{{\mathcal P}}
\newcommand{\sT}{{\tau}}

\newcommand{\Om}{{\Omega}}
\newcommand{\Sstyle}{\textstyle}
\DeclareMathOperator{\tr}{tr}
\DeclareMathOperator{\vol}{vol}

\newtheorem{theorem}{Theorem}

\makeatletter
\def\@sect#1#2#3#4#5#6[#7]#8{\ifnum #2>\c@secnumdepth
     \def\@svsec{}\else
     \refstepcounter{#1}\edef\@svsec{\csname the#1\endcsname.\hskip .75em }\fi
     \@tempskipa #5\relax
      \ifdim \@tempskipa>\z@
        \begingroup #6\relax
          \@hangfrom{\hskip #3\relax\@svsec}{\interlinepenalty \@M #8\par}%
        \endgroup
       \csname #1mark\endcsname{#7}\addcontentsline
         {toc}{#1}{\ifnum #2>\c@secnumdepth \else
                      \protect\numberline{\csname the#1\endcsname}\fi
                    #7}\else
        \def\@svsechd{#6\hskip #3\@svsec #8\csname #1mark\endcsname
                      {#7}\addcontentsline
                           {toc}{#1}{\ifnum #2>\c@secnumdepth \else
                             \protect\numberline{\csname the#1\endcsname}\fi
                       #7}}\fi
     \@xsect{#5}}
\def\@begintheorem#1#2{\it \trivlist \item[\hskip \labelsep{\bf #1\ #2.}]}
\makeatother

\begin{document}

\begin{center}
{\Large\bf Generalizations of Sch\"{o}bi's Tetrahedral Dissection }\\

\vspace*{+.5in}

N. J. A. Sloane \\
AT\&T Shannon Labs \\
180 Park Ave., Florham Park, NJ 07932-0971\\ [+.25in]

Vinay A. Vaishampayan \\
AT\&T Shannon Labs \\
180 Park Ave., Florham Park, NJ 07932-0971\\ [+.25in]

Email: njas@research.att.com, vinay@research.att.com
 
\vspace*{+.25in}

October 19, 2007; revised November 13, 2007

\vspace{+.5in}
{\bf Abstract}
\vspace*{+.2in}
\end{center}

%%%%%%%%% ABSTRACT

Let $v_1, \ldots, v_n$ be unit vectors in $\RR^n$ such that
$v_i \cdot v_j = -w$ for $i \ne j$ where
$-1 <w < \frac{1}{n-1}$.
The points $\sum_{i=1}^{n} \lambda_i v_i$
($1 \ge \lambda_1 \ge \cdots \ge \lambda_n \ge 0$)
form a ``Hill-simplex of the first type'',
denoted by $\sQ _n(w)$.
It was shown by Hadwiger in 1951 that $\sQ _n(w)$
is equidissectable with a cube.
In 1985, Sch\"{o}bi gave a three-piece dissection of 
$\sQ _3(w)$ into a triangular 
prism $c \sQ _2(\frac{1}{2}) \times I$,
where $I$ denotes an interval and $c = \sqrt{2(w+1)/3}$.
The present paper generalizes Sch\"{o}bi's dissection to
an $n$-piece dissection of $\sQ _n(w)$
into a prism $c \sQ _{n-1}(\frac{1}{n-1}) \times I$,
where $c = \sqrt{(n-1)(w+1)/n}$.
Iterating this process
leads to a dissection of $\sQ _n(w)$ into an $n$-dimensional
rectangular parallelepiped (or ``brick'') using at most $n!$ pieces.
The complexity of computing the map from $\sQ _n(w)$
to the brick is $O(n^2)$.
A second generalization of Sch\"{o}bi's dissection is
given which applies specifically in $\RR^4$.
The results have applications to source coding and to
constant-weight binary codes.

\vspace{0.8\baselineskip}
Keywords: dissections, Hill tetrahedra, Sch\"{o}bi, polytopes, 
Voronoi cell, source coding, constant-weight codes

\vspace{0.8\baselineskip}
2000 Mathematics Subject Classification: 52B45 (94A29, 94B60)

%%%% SECTION 1

\setlength{\baselineskip}{1.5\baselineskip}
\section{Introduction}\label{Sec1}

We define $\sQ _n(w)$ (where $n \ge 1$ and $-1 <w < \frac{1}{n-1}$)
as in the Abstract, and let 
$\sO _n := \sQ _n(0)$,
$\sP _n := \sQ _n(\frac{1}{n})$.
Hadwiger~\cite{Hadw1951} showed in 1951  
(see also Hertel~\cite{Hert2003})
that $\sQ _n(w)$ is equidissectable with a cube for all $n$.
His proof is indirect and not constructive.
The simplex $\sO_n$ is especially interesting:
it has vertices
\beql{EqH1}
000\ldots00,~
100\ldots00,~
110\ldots00,~
111\ldots00,~
\ldots,~
111\ldots10,~
111\ldots11  \,,
\eeq
and is an {\em orthoscheme} in Coxeter's terminology~\cite{Coxe1973}.
Because of applications to encoding
and decoding constant-weight codes~\cite{TVS2007},
we are interested in algorithms that carry out the
dissection of $\sO_n$ in an efficient manner.
In fact our question is slightly easier than the classical
problem, because we only need to decompose
$\sO_n$ into pieces which can be
reassembled to form a rectangular parallelepiped
(or $n$-dimensional ``brick''), not necessarily a cube\footnote{For the
problems of dissecting a rectangle into a square and
a three-dimensional rectangular parallelepiped into a cube see 
Boltianskii~\cite[p.~52]{Bolt1978},
Cohn~\cite{Cohn1974},
Frederickson~\cite[Page~236]{Fred1997}.}.

For the case $n=3$, Hill~\cite{Hill1895} had already shown in
1895 that the tetrahedra $\sQ_3(w)$ are equidissectable with a cube.
It appears that that the first explicit
dissection of $\sO_3$ into a cube was given by Sydler~\cite{Sydl1956}
in 1956.
Sydler shows that $\sO_3$ may be cut into four pieces
which can be reassembled to form a prism with base
an isosceles right triangle.
One further cut then gives a brick.
Sydler's dissection can be seen in a number of references
(Boltianskii~\cite[p.~99]{Bolt1978},
Cromwell~\cite[p.~47]{Crom1997},
Frederickson~\cite[Fig.~20.4]{Fred1997},
Sydler~\cite{Sydl1956},
Wells~\cite[p.~251]{Well1991})
and we will not reproduce it here.
Some of these references incorrectly attribute
Sydler's dissection to Hill.

In our earlier paper~\cite{TVS2007},
we gave a dissection of $\sO_n$ to a prism
$\sO_{n-1} \times I$ for all $n$
%$\sO_{n-1} \times I_{1/n}$ (where $I_{1/n}$ 
%denotes an interval of length $1/n$) for all $n$,
that requires $(n^2-n+2)/2$ pieces.
In three dimensions this uses four pieces,
the same number as Sydler's, but is somewhat simpler 
than Sydler's in that all our cuts are made along
planes perpendicular to coordinate axes.
By iterating this construction we eventually
obtain a dissection of $\sO_n$ into an $n$-dimensional brick.
The total number of pieces in the overall dissection
is large (roughly $(n!)^2/2^n$),
but the complexity of computing the coordinates
of a point in the final brick, given a initial point in $\sO_n$, is only
$O(n^2)$.

In 1985, Sch\"obi~\cite{Scho1985}\footnote{According
to Frederickson~\cite[Page 234]{Fred1997},
this construction was independently found by Anton Hanegraaf, unpublished.} gave
a dissection of $\sQ_3(w)$ (where $-1<w<\frac{1}{2}$)
into a prism with base an equilateral triangle that uses only three
pieces  (see Figs. \ref{Fig4new}, \ref{Fig5new} below,
also Frederickson~\cite[Fig.~20.5]{Fred1997}).
There is a way to cut $\sQ_n(w)$ for any $n$ into $n$ pieces
that is a natural generalization of Sch\"obi's dissection,
but for a long time we were convinced that
already for $n=4$ these pieces could not be reassembled to form 
a prism $P \times I$ for any $(n-1)$-dimensional polytope $P$.
In fact, we were wrong, and the main goal of this paper is to
use the ``Two Tile Theorem'' (Theorem \ref{Th1}) to generalize
Sch\"obi's dissection to all dimensions. 
%prove with the help of the ``Two Tile Theorem''
%that Sch\"obi's dissection can be generalized to
%all dimensions. 
We will show in Theorem \ref{Th2} that 
$\sQ_n(w)$ can be cut into $n$ pieces
that can be reassembled to form a prism 
$c \sP_{n-1} \times I_{\ell}$,
where
$c = \sqrt{ (n-1)(w+1)/n }$ and
$\ell = \sqrt{ (1-w(n-1))/n }$.
%$$
%\sqrt{ \frac{ (n-1)(w+1) }{n} } \sP_{n-1} \times I_{\ell} \,,
%$$
%where $ \ell := \sqrt{ \frac{ 1-w(n-1) }{n} } $.
The cross-section is always  proportional to
$\sP_{n-1} = \sQ_{n-1}(\frac{1}{n-1})$, independently of $w$.

By iterating this dissection we eventually decompose
$\sQ_n(w)$ (and in particular $\sO_n$) into a brick. 
The total number of pieces is at most $n!$
and the complexity of computing the map from $\sQ _n(w)$
to the brick is $O(n^2)$ (Theorem \ref{Th3}).
Although this is the same order of complexity as the algorithm 
given in our earlier paper \cite{TVS2007}, the present
algorithm is simpler and the number of
pieces is much smaller.

The recreational literature on dissections consists mostly 
of {\em ad hoc} constructions, although there are a few
general techniques, which can be found in the books of
Lindgren~\cite{Lind1964} and Frederickson~\cite{Fred1997}, \cite{Fred2002}.
The construction we have found the most useful is
based on group theory.
We call it the ``Two Tile Theorem'',
and give our version of it in Section \ref{Sec2},
together with several examples.
In Section \ref{Sec3} we state and prove the main theorem,
and then in Section \ref{Sec4} we study the overall algorithm for
dissecting $\sO_n$ into a brick.

Before finding the general dissection mentioned above, 
we found a different generalization of Sch\"obi's dissection
which applies specifically to the $4$-dimensional case.
This is described in Section \ref{Sec5}.
It is of interest because it is partially
(and in a loose sense) a ``hinged''
dissection (cf. Frederickson~\cite{Fred2002}).
After two cuts have been made, the first two motions
each leave a two-dimensional face fixed.
We then make a third cut, giving a total of six pieces
which can reassembled to give a prism $c \sP_3 \times I$.
This construction is also of interest because
it is symmetrical, and
it is the only {\em ad hoc} dissection we know of in four dimensions
(the dissections found by Paterson~\cite{Pate1996}
are all based on a version of the Two Tile Theorem).

A note about applications.
If we have a dissection of a polytope $P$ into a brick 
$I_{\ell_1} \times I_{\ell_2} \times  \cdots \times I_{\ell_n}$,
then we have a natural way to encode the points of $P$
into $n$-tuples of real numbers. This bijection provides a useful
parameterization of the points of $P$.
It may be used for source coding, if we
have a source that produces points uniformly
distributed over $P$ 
(for example, $P$ might be the Voronoi cell of a lattice).
Conversely, the bijection may be used in simulation,
when we wish to synthesize a uniform distribution
of points from $P$.
%Suppose we have a source that produces points
%that are uniformly distributed over a certain polytope $P$
%(for example, the Voronoi cell of a lattice).
%If there is a dissection of $P$ into a brick
%$I_{\ell_1} \times I_{\ell_2} \times  \cdots \times I_{\ell_n}$,
%then we have a natural way to encode the points of $P$
%into $n$-tuples of real numbers. This bijection provides a useful
%parameterization of the points of $P$.
For the application
to constant-weight codes we refer the
reader to \cite{TVS2007}.

\vspace*{+.1in}

\noindent{\bf Notation.}
A polytope in $\RR^n$ is a union of a finite 
number of finite $n$-dimensional simplices. It need be neither
convex nor connected.
Let $P, P_1, \ldots, P_k$ be polytopes in $\RR^n$.
By $P=P_1 + \cdots +  P_k$ we mean
that the interiors of $P_1, \ldots, P_k$
are pairwise disjoint and $P=P_1 \cup \ldots \cup P_k$.
Let $\Gamma$ be a group of isometries of $\RR^n$.
Two polytopes $P$, $Q$ in $\RR^n$ are said to be
\emph{$\Gamma$-equidissectable}
%, written $P \underset{\Gamma}{\thicksim} Q$,
if there are polytopes $P_1, \ldots, P_k$, $Q_1, \ldots, Q_k$
for some integer $k \ge 1$ such that
$P=P_1 + \ldots + P_k$,
$Q=Q_1 + \ldots + Q_k$, and
$P_1^{g_1} = Q_1, \ldots, P_k^{g_k} = Q_k$
for appropriate elements $g_1, \ldots, g_k \in \Gamma$.
In case $\Gamma$ is the full isometry group of $\RR^n$
we write $P~{\thicksim}~Q$ and
say that $P$ and $Q$ are \emph{equidissectable}. 
Isometries may involve reflections:
we do not insist that the dissections can be carried out
using only transformations of determinant $+1$.
$I_{\ell}$ denotes an interval of length $\ell$,
$I$ is a finite interval of unspecified length,
and $I_\infty = \RR^1$.

For background information about dissections and Hilbert's third problem,
and any undefined terms, we refer the reader to
the excellent surveys by
Boltianskii~\cite{Bolt1978},
Dupont~\cite{Dupo2001},
Frederickson~\cite{Fred1997}, \cite{Fred2002},
Lindgren~\cite{Lind1964},
McMullen~\cite{McMu1993},
McMullen and Schneider~\cite{McSc1983},
Sah~\cite{Sah1979}
and
Yandell~\cite{Yand2002}.

%%%% SECTION 2
\section{The ``Two Tile Theorem''}\label{Sec2}

Let $A \subset \RR^n$ be a polytope, $\Gamma$ a group
of isometries of $\RR^n$ and $\Om$ a subset of $\RR^n$.
If the images of $A$ under the action
of $\Gamma$ have disjoint interiors,
and $\Om = \cup_{g \in \Gamma} A^g$, we
say that $A$ is a \emph{$\Gamma$-tile} for $\Om$.
This implies that $\Gamma$ is discontinuous and 
fixed-point-free.

Versions of the following theorem---although
not the exact version that we need---have
been given by 
Aguil\'o, Fiol and Fiol~\cite[Lemma~2.2]{AFF2000},
M\"uller~\cite[Theorem~3]{Mull1988} and
Paterson~\cite{Pate1996}.
It is a more precise version of the technique of
%It is in the same spirit as Frederickson's
``superposing tesselations'' used by
Macaulay \cite{Maca1914}, \cite{Maca1919},
Lindgren \cite[Chap. 2]{Lind1964}
and
Frederickson \cite[p. 29]{Fred1997}, \cite[Chap. 3]{Fred2002}.

\begin{theorem}\label{Th1}
If for some set $\Om \subset \RR^n$ and
some group $\Gamma$ of isometries of $\RR^n$,
two $n$-dimensional polytopes $A$ and $B$ are both $\Gamma$-tiles for $\Om$,
then $A$ and $B$ are $\Gamma$-equidissectable.
\end{theorem}

\vspace*{+.1in}
\noindent{\bf Proof.}
We have
$$
A = A \cap \Om = A \cap \bigcup_{g \in \Gamma} B^g
= \bigcup_{g \in \Gamma} A \cap B^g \,,
$$
where only finitely many of the intersections $A \cap B^g$
are nonempty.
The set of nonempty pieces $\{ A \cap B^g \mid g \in \Gamma \}$
therefore gives a dissection of $A$,
and by symmetry
the set of nonempty pieces $\{ A^g \cap B \mid g \in \Gamma \}$
gives a dissection of $B$.
But $(A \cap B^g)^{g^{-1}} = A^{g^{-1}} \cap B$,
so the two sets of pieces are the same modulo 
isometries in $\Gamma$.~~~$\bsq$

We give four examples; the main application will
be given in the next section.

\vspace*{+.1in}

% This is Fig 1
\begin{figure}[htb]
\centerline{\includegraphics[width=5cm]{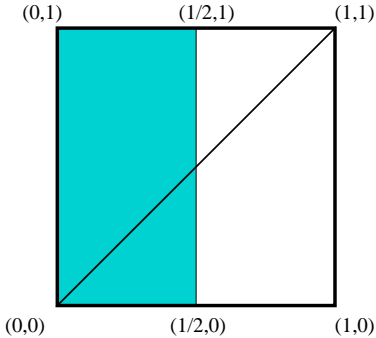}}
\caption{Illustrating the Two Tile Theorem: $A$ is the
triangle $(0,0), (1,0), (1,1)$, $B$ (shaded) is the rectangle
$(0,0), (\frac{1}{2},0),(\frac{1}{2},1),(0,1)$, $\Omega$ is the square
$(0,0),(1,0),(1,1),(0,1)$ and $\Gamma$ is 
generated by $\phi~:~(x,y) \mapsto (1-x,1-y)$.}
\label{trsq.fig}
\end{figure}

\noindent{\bf Example 1.}
Let $A=\sO_2$, the right triangle with vertices $(0,0), (1,0),
(1,1)$, let $B$ be the rectangle with vertices
$(0,0), (\frac{1}{2},0),(\frac{1}{2},1),(0,1)$ 
and let $\phi$ be the map $(x,y) \mapsto (1-x,1-y)$. Let
$\Gamma$ be the group of order $2$ generated by $\phi$ and let $\Omega$
be the square $(0,0),(1,0),(1,1), (0,1)$. Then $A$ and $B$ are
both $\Gamma$-tiles for $\Omega$. It follows
from Theorem \ref{Th1} that $A$ and $B$ are
equidissectable (see
%with pieces as shown in
Fig.~\ref{trsq.fig}).
%Observe that $\Omega$ is an invariant set for the group generated
%by $\phi$.
Alternatively, we could take the origin to be at the center
of the square, and then the theorem applies with
$\phi := (x,y) \mapsto (-x,-y)$.

%% This is Fig 2
\begin{figure}[htb]
%\begin{center}
\centerline{\includegraphics[width=9.5cm]{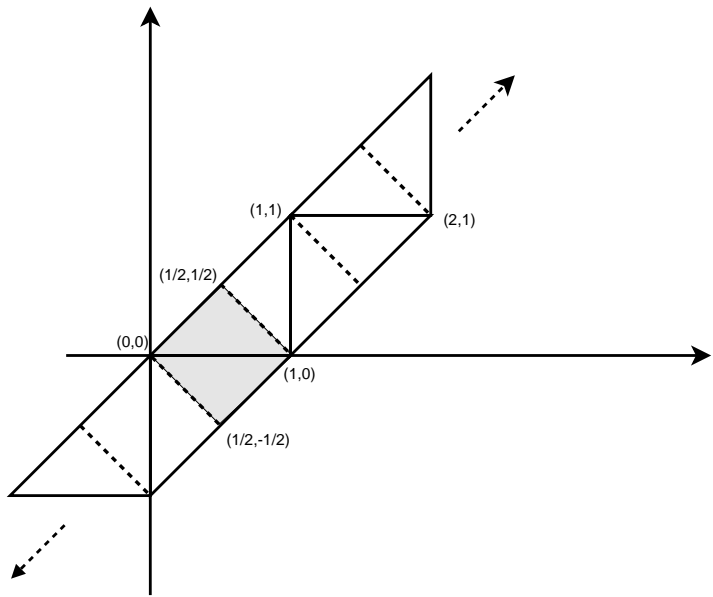}}
%\centerline{\includegraphics[width=7.5cm]{a110312_3v.eps}}
%\end{center}
\caption{ Another illustration of the Two Tile Theorem: $A$ is
the triangle $(0,0), (1,0), (1,1)$,
$B$ (shaded) is the square $(0,0),
(\frac{1}{2},-\frac{1}{2}),
(1,0),
(\frac{1}{2},\frac{1}{2})$
and $\Om$ is the strip $x \ge y \ge x-1$. }
\label{Fig2a}
\end{figure}

%\vspace*{+.1in}

\noindent{\bf Example 2.}
Again we take $A = \sO_2$ to be the right triangle
with vertices $(0,0), (1,0), (1,1)$,
but now we take $\phi$ to be the map $(x,y) \mapsto (y+1,x)$.
Note that $\phi$ involves a reflection.
As mentioned in \S\ref{Sec1}, this is permitted by 
our dissection rules.
Let $\Gamma$ be the infinite cyclic group generated by $\phi$,
and let $\Om$ be the infinite strip
defined by $x \ge y \ge x-1$.
Then $A$ is a $\Gamma$-tile for $\Om$ (see Fig. \ref{Fig2a}).
For $B$, the second tile, we take the square
with vertices $(0,0), (\frac{1}{2},-\frac{1}{2}),
(1,0), (\frac{1}{2},\frac{1}{2})$. This is also a $\Gamma$-tile for $\Om$,
and so $A$ and $B$ are equidissectable.
The two nonempty pieces are the triangles $A \cap B$
and $A \cap B^{\phi}$. The latter is mapped
by $\phi^{-1}$ to the triangle with vertices
$(0,0), (\frac{1}{2},-\frac{1}{2}), (1,0)$.
This is a special case of the dissection given in Theorem \ref{Th2}.
Of course in this case the dissection could
also have been accomplished without using reflections.

%% This is Fig 4
\begin{figure}[htb]
\begin{center}
\includegraphics[width=7.5cm]{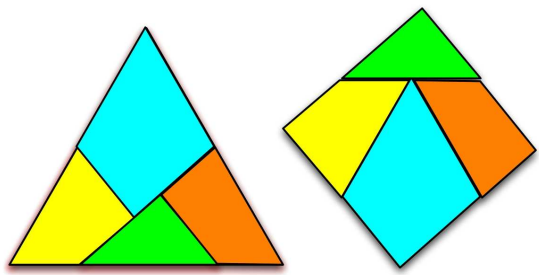}
\end{center}
\caption{ Four-piece dissection of an
equilateral triangle to a square, usually attributed 
to Dudeney (1902)}
\label{Fig2b}
\end{figure}

%\vspace*{+.1in}
\noindent{\bf Example 3.}
One of the most elegant of all dissections is the well-known
four-piece dissection of an 
equilateral triangle to a square, shown in Fig. \ref{Fig2b}.
This was published in 1902 by Dudeney, although
Frederickson~\cite[Page~136]{Fred1997} suggests that
he may not have been the original discoverer.
This dissection can be found in many references
(for example, Coffin~\cite[Chap.~1]{Coff1991},
Eves~\cite[\S5.5.1]{Eves1966},
%Frederickson~\cite[p.~2]{Fred2002},
Wells~\cite[p.~61]{Well1991}).
Gardner \cite[Chap.~3]{Gard1961} gives a proof by elementary geometry.
The usual construction of this dissection, however, is 
by superimposing two strips, a technique that 
Lindgren calls a $TT$-dissection
(Akiyama and Nakamura \cite{AkNa1998},
Frederickson \cite[Chaps.~11,~12]{Fred1997}, \cite[Chap.~3]{Fred2002},
Lindgren \cite[Fig.~5.2]{Lind1964}).
The literature on dissections does not appear to contain
a precise statement of conditions which guarantee that 
this construction produces a dissection.
Such a theorem can be obtained as a corollary of the
Two Tile Theorem, and will be published elsewhere, together
with rigorous versions of other strip dissections.
Both Gardner and Eves mention that L.~V.~Lyons
extended Dudeney's dissection to 
cut the whole plane into a ``mosaic of interlocking
squares and equilateral triangles,'' and Eves shows this ``mosaic''
in his Fig.~5.5b
(Fig.~\ref{Fig2c} below shows essentially the same figure, with the addition
of labels for certain points).
We will use this ``mosaic,'' which is really a double tiling
of the plane, to give an alternative proof that the
dissection is correct from the Two Tile Theorem.
Following Lyons, we first use the dissection to construct the 
double tiling. We then ignore how the double
tiling was obtained, and apply the Two Tile Theorem
to give an immediate certificate of proof for Dudeney's dissection.
The double tiling also has some interesting properties
that are not apparent from Eves's figure,
and do not seem to have been mentioned before
in the literature.
%However, {\em once we are told that the dissection exists},
%we can give an alternative derivation,
%and hence a certificate of proof, for it,
%using the Two Tile Theorem.
%The resulting double tiling of the plane
%has some interesting properties.

%% This is Fig 5
\begin{figure}[htb]
%\begin{center}
\centerline{\includegraphics[width=14.5cm]{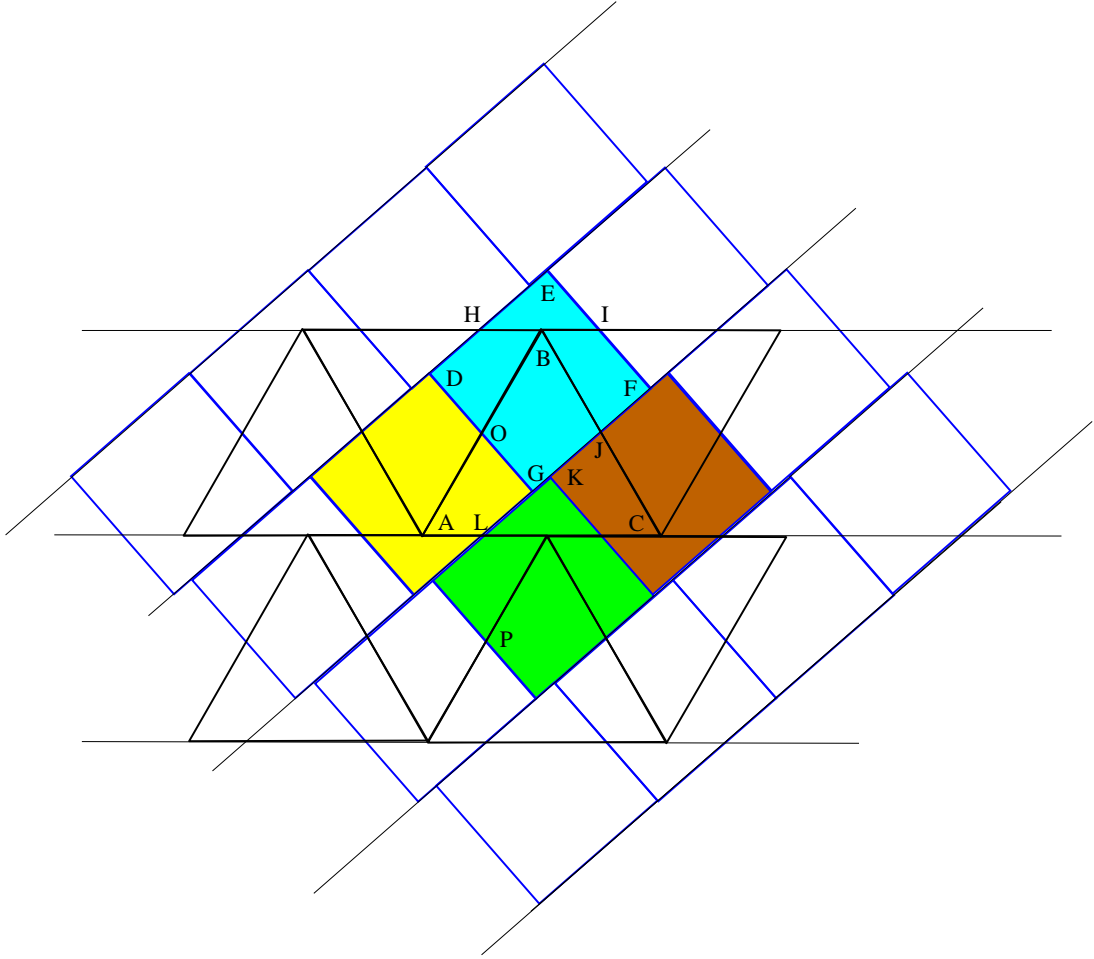}}
%\end{center}
\caption{ 
Lyons's ``mosaic,'' a double tiling
of the plane by triangles and squares.
%The dissection of Fig. \ref{Fig2b} derived from
%the Two Tile Theorem.
}
\label{Fig2c}
\end{figure}

Let $\Omega = \RR^2$, and
take the first tile to be an equilateral triangle
with edge length $1$, area $c_1 := \frac{\sqrt{3}}{4}$
and vertices 
$A := (-1/4, -c_1)$, $B := -A$ and
%$B := ( 1/4,  c_1)$ and
$C := ( 3/4, -c_1)$ (see Fig. \ref{Fig2c}),
with the origin $O$ at the midpoint of $AB$.
The second tile is a square with edge length $c_2 := \sqrt{c_1}$.
The existence of the dissection imposes many constraints, such as
$|JB| = |JC| = |HI| = 1/2$,
$|OD| = |OG| = c_2/2$,
$2|LG| + |GK| = 2|JK| +|GK| = c_2$, etc.,
and after some calculation we find that the
square should have vertices
$D := (     -c_1/2,  c_3/2      )$,
$E := (c_3 - c_1/2,  c_3/2 + c_1)$,
$F := (c_3 + c_1/2, -c_3/2 + c_1)$ and
$G := -D$,
where $c_3 = c_2 \sqrt{1-c_1}$.

We now construct a strip of squares that replicates the square
$DEFG$ in the southwest/northeast direction,
and a strip of triangles replicating $ABC$ (with alternate triangles
inverted) in the horizontal direction.
In order to fill the plane with copies of these strips, we 
must determine the offset of one strip of squares with respect to the
next strip of squares, and of one strip
of triangles with respect to the
next strip of triangles.
This implies the further constraints that
$P-O = L-H = K-E$, etc.,
and in particular that
$P$ should be the point $(1-2c_3, -2c_1)$.
Other significant points are
$H := -L := (c_3-1/2,c_1)$, $I := (c_3,c_1)$,
$J := (1/2,0)$, $K := (1-c_3-c_1/2, c_3/2-c_1)$.
The angle $CLG$ is $\arctan(c_1/c_3) = 41.15\ldots$ degrees.

We now have the desired double tiling of the plane.
Both the triangle $ABC$ and the square $DEFG$
are $\Gamma$-tiles for the whole plane,
where $\Gamma$ is the group
(of type $p2$ in the classical notation, 
or type $2222$ in the orbifold notation)
generated by translation by $(1,0)$,
translation by $OP = (1-2c_3, -2c_1) = (0.009015\ldots, -0.8660\ldots)$,
and multiplication by $-1$.
We now ignore how this tiling was found, and deduce from
Theorem \ref{Th1} that Dudeney's dissection exists.
The four pieces are $OBJG$, $ODHB$, $HEI$ and $BIFJ$.

It is interesting that the horizontal strips
of triangles do not line up exactly: each strip is shifted to
the left of the one below it by
$1-2c_3 = 0.009015\ldots$.
The group is correspondingly more complicated than one might have
expected from looking at Fig. \ref{Fig2c},
since the second generator for the group is not {\em quite} translation
by $(0, -\sqrt{3}/2)$!
There is an associated lattice, generated by the vectors $OH$ and $OJ$,
and containing the points $I$ and $L$,
which is {\em nearly} rectangular, the angle between
the generators being $89.04\ldots$ degrees.

Incidentally, although Lindgren
\cite[p.~25]{Lind1964}
refers to this dissection as ``minimal'', we have never seen
a proof that a three-piece dissection of an equilateral
triangle to a square is impossible.
This appears to be an open question.

\vspace*{+.1in}
\noindent{\bf Example 4.}
It is easy to show by induction that any lattice
$\Lambda$ in $\RR^n$ has a brick-shaped fundamental
region.  The theorem then provides a dissection of the Voronoi
cell of $\Lambda$ into a brick.
For example, the Voronoi cell of the root lattice $D_n$ is described
in \cite[Chap. 21]{SPLAG}. By applying the theorem,
we obtain a dissection of the Voronoi cell into a brick that uses $2n$ pieces.
For $n=3$ this gives the well-known six-piece dissection of a rhombic
dodecahedron into a $2 \times 1 \times 1$ brick
(cf. \cite[pp.~18,~242]{Fred1997}).
%is stated on p. 242 that seven pieces are needed, but two of these
%pieces can be combined.}).
%Details will be published elsewhere.

%%%% SECTION 3
\section{The main theorem}\label{Sec3}

We begin by choosing a particular realization of the 
simplex $\sQ_n(w)$.
Define the following vectors in $\RR^n$:
\beql{EqVV}
v_1 := (a,b,b,\ldots,b), \,
v_2 := (b,a,b,\ldots,b), \,
v_3 := (b,b,a,\ldots,b), \,
\ldots, \,
v_n := (b,b,b,\ldots,a),
\eeq
where
\begin{eqnarray}
b & := &  (\sqrt{1-w(n-1)}-\sqrt{1+w})/n \,, \nonumber  \\
a & := & b+\sqrt{1+w} \,. 
\label{Eqab}
\end{eqnarray}
Then $v_i \cdot v_i = 1$, $v_i \cdot v_j = -w$
for $i \ne j$, $i,j = 1, \ldots, n$.
We take the convex hull of the vectors
$0, v_1, v_1+v_2, \ldots, v_1+ \cdots + v_n$,
that is, the zero vector together with the rows of
\beql{Eqv1}
\left[ \begin{array}{ccccc}
a &b &b & \ldots &b \\
a+b &a+b &2b & \ldots &2b \\
a+2b &a+2b &a+2b & \ldots &3b \\
\ldots &\ldots &\ldots & \ldots &\ldots \\
a+(n-1)b &a+(n-1)b &a+(n-1)b & \ldots &a+(n-1)b \end{array} \right] \,,
\eeq
to be our standard version of $\sQ_n(w)$.
This simplex has volume 
\beql{EqVol}
(1+w)^{(n-1)/2} (1-w(n-1))^{1/2} / n! \,.
\eeq
Setting $w=0, a=1, b=0$
%in \eqn{EqVV}
gives our standard version of $\sO_n$, as in \eqn{EqH1},
and setting $w=1/n, a=n^{-3/2}(1+(n-1)\sqrt{n+1}),
b=n^{-3/2}(1-\sqrt{n+1})$
gives our standard version of $\sP_n$.

Two other versions of  $\sP_n$ will also appear.
Let $p_i := 1/\sqrt{i(i+1)}$, 
%for $i=1,\ldots$,
and construct a $n \times n$ orthogonal matrix $M_n$
as follows. For $i=1,\ldots,n-1$ the $i$th
column of $M_n$ has entries
$p_i$ ($i$ times), $-ip_i$ (once) and $0$ ($n-i-1$ times),
and the entries in the last column are all $1/\sqrt{n}$. 
(The last column is in the 
$(1,1,\ldots,1)$
direction and the other columns are perpendicular to it.)
For example,
\renewcommand{\arraystretch}{1.25}
$$
M_3 := \left[ \begin{array}{ccc}
\frac{1}{\sqrt{2}} & \frac{1}{\sqrt{6}} & \frac{1}{\sqrt{3}} \\
-\frac{1}{\sqrt{2}} & \frac{1}{\sqrt{6}} & \frac{1}{\sqrt{3}} \\
0 & -\frac{2}{\sqrt{6}} & \frac{1}{\sqrt{3}} \end{array} \right] \,.
$$
\renewcommand{\arraystretch}{1}
The other two versions of  $\sP_n$ are: the convex hull
of the zero vector in $\RR^n$ together with the rows of
\beql{Eqv2}
\sqrt{\frac{n+1}{n}} \,
\left[ \begin{array}{ccccc}
p_1 & p_2 & p_3 & \ldots & p_n \\
0   & 2p_2 & 2p_3 & \ldots & 2p_n \\
\ldots &\ldots &\ldots & \ldots &\ldots \\
0 &0 &0 & \ldots & dp_n \end{array} \right] \,,
\eeq
and  the convex hull of
the zero vector in $\RR^{n+1}$ together with the rows of
\renewcommand{\arraystretch}{1.25}
\beql{Eqv3}
\sqrt{\frac{n+1}{n}} \,
\left[ \begin{array}{ccccc}
\frac{n}{n+1} & -\frac{1}{n+1} & -\frac{1}{n+1} & \ldots & -\frac{1}{n+1} \\
\frac{n-1}{n+1} & \frac{n-1}{n+1} & -\frac{2}{n+1} & \ldots & -\frac{2}{n+1} \\
\ldots &\ldots &\ldots & \ldots &\ldots \\
\frac{1}{n+1} & \frac{1}{n+1} & \frac{1}{n+1} & \ldots & -\frac{n}{n+1} 
\end{array} \right] \,.
\eeq
\renewcommand{\arraystretch}{1}
To see that both of these simplices
are congruent to the standard version of $\sP_n$,
note that multiplying \eqn{Eqv3} on the right by
$M_{n+1}$ produces \eqn{Eqv2} supplemented by a column of zeros,
and then multiplying \eqn{Eqv2} on
the right by $M_n^{\tr}$ (where tr denotes transpose)
produces the standard version.

\noindent{\bf Remark.}
If we ignore for the moment the scale factor in front
of \eqn{Eqv3}, we see that its rows are the coset
representatives for the root lattice $A_n$
in its dual $A_n^{\ast}$ \cite[p.~109]{SPLAG}.
In other words, the rows of \eqn{Eqv3}
contain one representative of each of the classes
of vertices of the Voronoi cell for $A_n$.
$\sP_2$ is an equilateral triangle and
$\sP_3$ is a ``Scottish tetrahedron''
in the terminology of Conway and Torquato~\cite{CoTo2006}.

We can now state our main theorem.

%%%%%% MAIN

\begin{theorem}\label{Th2}
The simplex $\sQ_n(w)$ is equidissectable 
with the prism $c \sP_{n-1} \times I_{\ell}$,
where \\
$c := \sqrt{ (n-1)(w+1)/n }$ and
$\ell := \sqrt{ (1-w(n-1))/n }$.
\end{theorem}

\vspace*{+.1in}
\noindent{\bf Proof.}
Let $\Om$ be the convex hull of the points 
$\{u_i \in \RR^n \mid i \in \ZZ \}$,
where $u_0:=(0,0,\ldots,0)$,
$u_i:=u_0^{\phi^i}$,
$\phi$ is the map
$$
\phi: (x_1,\ldots,x_n) \mapsto
(x_n+a, x_1+b, x_2+b, \ldots,x_{n-1}+b)
$$
and $a,b$ are as in \eqn{Eqab}
(see Table \ref{Tab1}).

\begin{table}[htb]
\caption{Points defining the infinite prism $\Om$.
The convex hull of any
$n+1$ successive rows is a copy of $\sQ_n(w)$.}
$$
\begin{array}{|ccccccc|}
\hline
\ldots & \ldots & \ldots & \ldots & \ldots & \ldots & \ldots \\
u_{-1} & = & -b & -b & -b & \ldots & -a \\
u_0 & = & 0 & 0 & 0 & \ldots & 0 \\
u_{1} & = & a & b & b & \ldots & b \\
u_{2} & = & a+b & a+b & 2b & \ldots & 2b \\
\ldots & \ldots & \ldots & \ldots & \ldots & \ldots & \ldots \\
u_{n} & = & a+(n-1)b & a+(n-1)b & a+(n-1)b & \ldots & a+(n-1)b \\
u_{n+1} & = & 2a+(n-1)b & a+nb & a+nb & \ldots & a+nb \\
\ldots & \ldots & \ldots & \ldots & \ldots & \ldots & \ldots  \\
\hline
\end{array}
$$
\label{Tab1}
\end{table}

We now argue in several easily verifiable steps.

\noindent{(i)}
For any $i \in \ZZ$, the convex hull of
$u_i, u_{i+1}, \ldots, u_{i+n}$ is a copy of 
$\sQ_n(w)$, $\sQ ^{(i)}$ (say),
with $\sQ ^{(i)} \subset \Om$ and
$(\sQ ^{(i)})^{\phi} = \sQ^{(i+1)}$.

\noindent{(ii)}
The simplices $\sQ ^{(i)}$ and $\sQ ^{(i+1)}$
share a common face, the convex hull of
$u_{i+1}, \ldots, u_{i+n}$, but have disjoint interiors.
More generally, for all $i \ne j$, $\sQ ^{(i)}$ and $\sQ ^{(j)}$
have disjoint interiors.

\noindent{(iii)}
The points of $\Om$ satisfy
\beql{EqWall}
x_1 \ge x_2 \ge \cdots \ge x_n \ge x_1 - (a-b) \,.
\eeq
(This is true for $\sQ ^{(0)}$ and the property
is preserved by the action of $\phi$.)
The inequalities \eqn{EqWall} define an infinite prism
with axis in the 
$(1,1,\ldots,1)$
direction.
We will show that every point in the prism belongs to
$\Om$, so $\Om$ is in fact equal to this prism.

\noindent{(iv)}
The projection of $\sQ ^{(0)}$ onto the hyperplane
perpendicular to the $(1,1,\ldots,1)$ direction
is congruent to $c \sP_{n-1}$, where
$c := \sqrt{ (n-1)(w+1)/n }$.
(For multiplying \eqn{Eqv1} on the right
by $M_n$ gives a scaled copy of \eqn{Eqv2}.)
On the other hand, the intersection of the prism
defined by \eqn{EqWall} with the hyperplane $\sum_{i=1}^{n} x_i = 0$
consists of the points $(0,0,\ldots,0)$,
$ \sqrt{w+1} ((n-1)/n, -1/n, \ldots, -1/n)$,
$ \sqrt{w+1} ((n-2)/n, (n-2)/n, -2/n, \ldots, -2/n)$,
$\ldots$,
and---compare \eqn{Eqv3}---is also congruent to $c \sP_{n-1}$.
Since the projection and the intersection have
the same volume, it follows that every point in
the prism is also in $\Om$. (For consider a long
but finite segment of the prism. The total volume of
the copies of $\sQ_n(w)$ is determined by the projection,
and the volume of the prism is determined by the cross-section,
and these coincide.) We have therefore established that $\Om$ is
the infinite prism 
$$c \sP_{n-1} \times I_{\infty}$$
with walls given by \eqn{EqWall}.
Furthermore, $\sQ_n(w)$ is a $\Gamma$-tile for $\Om$,
where $\Gamma$ is the infinite cyclic group generated by $\phi$.

\noindent{(v)}
For a second tile, we take the prism
$$ B := c \sP_{n-1} \times I_{\ell} \,,$$
where 
$\ell := \sqrt{ (1-w(n-1))/n }$.
The length $\ell$ is chosen so that $B$ has the same volume as 
$\sQ_n(w)$ (see \eqn{EqVol}).
We take the base of $B$ to be the particular copy of $c \sP_{n-1}$
given by the intersection of $\Om$ with
the hyperplane $\sum_{i=1}^{n} x_i = 0$, as in (iv).
The top of $B$ is found by adding $\ell/\sqrt{n}$ to 
every component of the base vectors. 
To show that $B$ is also a $\Gamma$-tile for $\Om$,
we check that the image of the base of $B$
under $\phi$ coincides with the top of $B$. This is an easy
verification.

Since $\sQ_n(w)$ and $B$ are both $\Gamma$-tiles for $\Om$,
the desired result follows from Theorem \ref{Th1}.~~~$\bsq$

\noindent{\bf Remarks.}
(i) The prism $B$ consists of the portion
of the infinite prism $\Om$ bounded by the
hyperplanes $\sum x_i = 0$
and $\sum x_i = \sqrt{1-w(n-1)}$.
The ``apex'' of $\sQ_w(n)$ is the point
$(a+(n-1)b, a+(n-1)b, \ldots, a+(n-1)b)$,
which---since $a+(n-1)b = \sqrt{1-w(n-1)}$---lies on the hyperplane
$\sum x_i = n\sqrt{1-w(n-1)}$.
There are therefore $n$ pieces $\sQ_n(w) \cap B^{{\phi}^k}$
($k= 0,1,\ldots,n-1$)
in the dissection, obtained by cutting $\sQ_w(n)$
along the hyperplanes $\sum x_i = k\sqrt{1-w(n-1)}$
for $k=1,\ldots,n-1$.
To reassemble them to form $B$, we apply
$\phi^{-k}$ to the $k$th piece.

% This is the new Fig 4
\begin{figure}[htb]
\begin{center}
\centerline{\includegraphics[width=9.5cm]{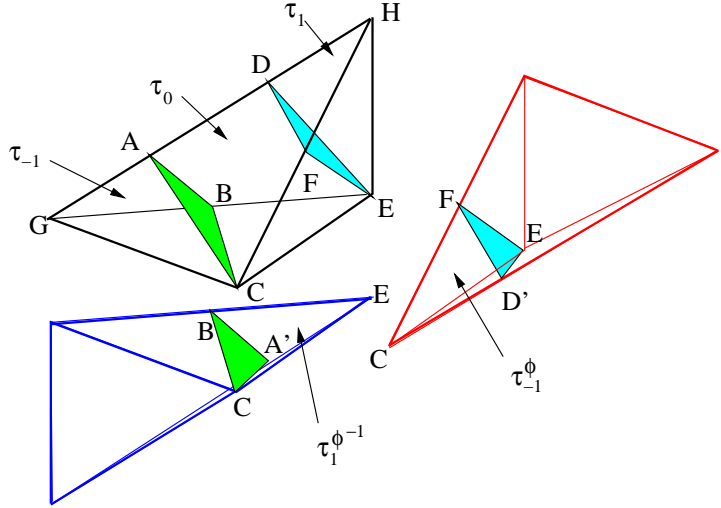}}
\end{center}
\caption{Exploded view showing three adjacent copies of $\sQ_3(w)$ and
their intersections with the two cutting planes.}
\label{Fig4new}
\end{figure}

\begin{figure}[htb]
\begin{center}
\centerline{\includegraphics[width=5.5cm]{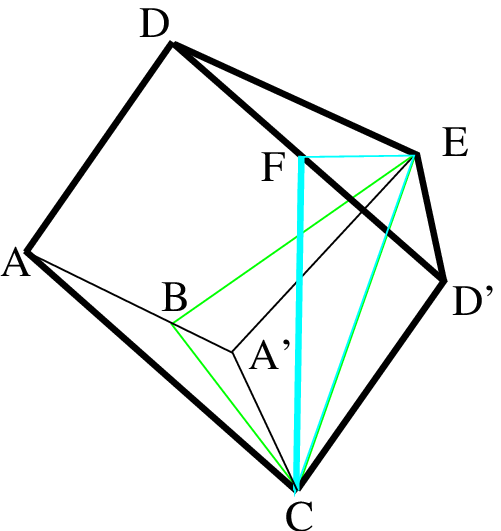}}
\end{center}
\caption{  Three-piece dissection of the triangular prism.
Our construction and Sch\"obi's use
the same three pieces but assemble them in a different way. In
this view the points $A'$ and $B$ are at the back of the figure. }
%\end{center}
%\label{Fig3a}
\label{Fig5new}
\end{figure}

(ii) The case $n=2$, $w=0$ of the theorem was
illustrated in Fig. \ref{Fig2a}.
In the case $n=3$, $-1~<w<~\frac{1}{2}$,
the three pieces are exactly the same as those in 
Sch\"obi's dissection \cite{Scho1985}.
However, it is interesting that we reassemble them
in a different way to form the same prism $c \sP_2 \times I_{\ell}$,
with $c := \sqrt{2(w+1)/3}$, $\ell := \sqrt{(1-2w)/3}$.
First we describe our dissection, which is illustrated in Fig. \ref{Fig4new}.
The figure shows an exploded view of three adjacent copies of $\sQ_3(w)$,
namely
$\sQ_3(w)^{\phi^{-1}}$ (the lower left tetrahedron), 
$\sQ_3(w)$ (the upper left tetrahedron) and
$\sQ_3(w)^{\phi}$ (the tetrahedron on the right),
and their intersections with the two cutting planes. 
The three pieces in the dissection can be seen in the
upper left tetrahedron $\sQ_3(w)$. They are
$ \sT_{-1} = \sQ_3(w) \cap B^{\phi^{-1}}$ (the piece $ABCG$ on the left 
of this tetrahedron),
$ \sT_0 = \sQ_3(w) \cap B$ (the central piece $ABCDEF$)
and $ \sT_1 = \sQ_3(w) \cap B^{\phi}$ (the piece $DEFH$ on the right).
In Fig. \ref{Fig4new} we can also see an exploded view of these three pieces
reassembled to form the triangular prism:
$\sT_1^{\phi^{-1}}$ is the right-hand piece $A'BCE$ of the lower figure and
$\sT_{-1}^{\phi}$ is the left-hand piece $CD'EF$ of the figure on the right.
The fully assembled prism is shown in Fig. \ref{Fig5new}:
the tetrahedron $A'BCE$ is $\sT_1^{\phi^{-1}}$ and
the tetrahedron $CD'EF$ is $\sT_{-1}^{\phi}$.

% This is the new Fig 5

On the other hand, Sch\"obi
reassembles the same pieces by 
rotating $ \sT_1$ about the edge $EF$ (which acts as a hinge),
sending $D$ to $D'$ and giving the tetrahedron $CD'EF$, and
rotating $ \sT_{-1}$ about the hinge $BC$,
sending $A$ to $A'$ and giving the tetrahedron $A'BCE$.
This is strictly different from our construction,
since $\phi$ has no fixed points.
The pieces are the same and the end result is the same,
but the two outer pieces $\sT_1$ and $\sT_{-1}$
have been interchanged!
% let us regard $\sT_1$ as the fixed piece.
% Then we would form the prism as
% \beql{Eq3D}
% \sT_0^{\phi} + \sT_1 + \sT_2^{\phi^{-1}} \,.
% \eeq
% On the other hand, Sch\"obi proceeds by
% {\em rotating} $\sT_0$ and $\sT_2$
% about certain fixed edges (which
% act as hinges) until they just touch $\sT_1$.
% This apparent paradox is explained by the fact that 
% when rotated about the hinge, $\sT_0$
% becomes $\sT_2^{\phi^{-1}}$,
% and similarly $\sT_2$ becomes $\sT_0^{\phi}$.
% So Sch\"obi's assembly of the prism is also
% given by \eqn{Eq3D},
% albeit by a different pair of motions!
 
%%% this is the figure we decided to drop
%%% This is Fig 3a
%%\begin{figure}[htb]
%%\begin{center}
%\centerline{\includegraphics[width=7.5cm]{TetraPrism3.eps}}
%%\centerline{\includegraphics[width=7.5cm]{Tetra3Exploded.eps}}
%%\end{center}
%%\caption{ 
%%Sch\"obi's three-piece dissection of $\sQ_3(w)$
%%to a triangular prism. Theorem \ref{Th2}
%%produces the same pieces but assembles them
%%in a different way. }
%\caption{The image of $\sQ_3(w)$ under the group G. Every unmarked tetrahedron
%can be obtained from one of the marked tetrahedra by a tranlation in the $(1,1,\ldots,1)$ direction.}
%\label{Fig3a}
%\end{figure}

(iii) By repeated application of Theorem \ref{Th2}
we can dissect $\sQ_n(w)$ into an $n$-dimensional
brick. Each of the $n$ pieces from the first
stage is cut into at most $n-1$ pieces at the second stage, and so on,
so the total number of pieces in the final
dissection is at most $n!$.
(It could be less, if a piece from one stage is not intersected
by all of the cutting planes at the next stage. 
It seems difficult to determine the exact number of pieces.) 

%%%% SECTION 4
\section{Dissecting $\sO_n$ into a brick}\label{Sec4}

In this section we discuss in more detail the
recursive dissection of $\sQ_n(w)$
into a rectangular parallelepiped
or ``brick'' in the case of greatest interest to
us, when we start with $\sO_n = \sQ_n(0)$.

From Theorem \ref{Th2} we have
\begin{eqnarray}
\sO_n & ~{\thicksim}~ & \sqrt{\frac{n-1}{n}} \,\sP_{n-1} \times I_{\frac{1}{\sqrt{n}}}  \,, \nonumber  \\
\sP_n & ~{\thicksim}~ &  \frac{\sqrt{n^2-1}}{n} \sP_{n-1} \times I_{\frac{1}{n}} \,, 
\label{Eq84}
\end{eqnarray}
and so (since $\sP_1 = I_1$)
\beql{Eq85}
\sO_n ~{\thicksim}~
\frac{1}{2} I_1 \times I_{\Sstyle {p_2} } \times I_{\Sstyle {p_3} }
\times \cdots \times
I_{\Sstyle {p_{n-1}} } \times I_{\frac{1}{\sqrt{n}}} \,.
\eeq
The right-hand side of \eqn{Eq85} is our final brick; 
we will denote it by $\Pi$. Note that 
$\vol (\sO_n) = \vol (\Pi) = 1/n!$.

Let $\Theta$ denote the map from $\sO_n$ to $\Pi$
associated with the dissection \eqn{Eq85}.
We will show that given $x := (x_1, \ldots, x_n) \in \sO_n$,
$(y_1, \ldots, y_n) := \Theta(x) \in \Pi$
can be computed in $O(n^2)$ steps.

The algorithm for computing $\Theta$ breaks up
naturally into two parts.
The first step involves
dissecting $\sO_n$ into $n$ pieces
and reassembling them to form the prism
$$
B := \sqrt{\frac{n-1}{n}} \,\sP_{n-1} \times I_{\frac{1}{\sqrt{n}}} \,.
$$
All later steps start with a point in $\lambda_k \sP_{k}$
for $k = n-1, n-2, \ldots, 2$ and certain constants $\lambda_k$,
and produce a point in $\lambda_{k-1} \sP_{k-1} \times I$.

For the first step we must determine which 
of the pieces $\sO_n \cap B^{\phi_1^r}$
($r= 0,1,\ldots,n-1$)
$x$ belongs to, where $\phi_1$ is the map
$(x_1,\ldots,x_n) \mapsto (x_n+1, x_1, x_2, \ldots,x_{n-1})$.
This is given by $r := \lfloor \sum_{i=1}^{n} x_i \rfloor$,
and then mapping $x$ to
$x' := x^{\phi_1^{-r}}$ corresponds to reassembling the pieces
to form $B$.  However, $x'$ is expressed in terms of 
the original coordinates for $\sO_n$ and we must multiply it by $M_n$
to get coordinates perpendicular to the $(1,1,\ldots,1)$
direction, getting $x'' := (x_1'', \ldots, x_{n-1}'', y_n) = x' M_n$.
The final component of $x''$ is the projection of $x'$ 
in the $(1,1,\ldots,1)$ direction.
The other components of $x''$,
$(x_1'', \ldots, x_{n-1}'')$ define a point in 
$\sqrt{\frac{n-1}{n}} \,\sP_{n-1}$, but expressed in coordinates of the 
form shown in \eqn{Eqv2}, and before we proceed to the next stage,
we must reexpress this in the standard coordinates for
$\sqrt{\frac{n-1}{n}} \,\sP_{n-1}$, which we do by
multiplying it by
$M_{n-1}^{\tr}$ (see the beginning of \S\ref{Sec3}),
getting $x'''$.

The following pair of observations shorten these calculations.
First, $y_n$ can be computed 
directly once we know $r$, since each application of $\phi_1^{-1}$
subtracts $1$ from the sum of the coordinates.
If $s := \sum_{i=1}^{n} x_i$, then $r := \lfloor s \rfloor$
and $y_n = (s-r)/\sqrt{n}$.
Second,  the product of $M_n$-with-its-last-column-deleted 
and $M_{n-1}^{\tr}$ is the $n \times (n-1)$ matrix
\renewcommand{\arraystretch}{1.25}
$$ 
N_n := \left[ \begin{array}{cccc}
1-p_n & -p_n & \ldots & -p_n \\
-p_n & 1-p_n & \ldots & -p_n \\
\ldots & \ldots & \ldots & \ldots \\
-p_n & -p_n & \ldots & 1-p_n \\
-\frac{1}{\sqrt{n}} & -\frac{1}{\sqrt{n}} & \ldots & -\frac{1}{\sqrt{n}}
\end{array} \right] \,.
$$
\renewcommand{\arraystretch}{1}
Multiplication by $N_n$ requires only $O(n)$ steps.

The first stage in the computation of $\Theta$
can therefore be summarized as follows:

\begin{quotation}
\noindent
Step ${\rm A}$. Given $x:=(x_1, \ldots, x_n) \in \sO_n$. 
Let $s := \sum_{i=1}^{n} x_i$, $r := \lfloor s \rfloor$. \\
Compute $x' := x^{\phi_1^{-r}}$. \\
Pass $x''' := x' N_n$ to the next stage, and output 
$y_n := (\mbox{fractional~part~of~} s)/\sqrt{n}$.
\end{quotation}

In all the remaining steps we start with a point $x$
in $\lambda_k \sP_k$ for some constant $\lambda_k$,
where $k=n-1$, $n-2, \ldots,2$. Instead of $\phi_1$ we use
the map
$\phi_2 : (x_1,\ldots,x_k) \mapsto (x_k+a, x_1+b, x_2+b, \ldots,x_{k-1}+b)$,
where 
$a=k^{-3/2}(1+(k-1)\sqrt{k+1}),
b=k^{-3/2}(1-\sqrt{k+1})$.
Each application of $\phi_2^{-1}$ subtracts $1/\sqrt{k}$ from
the sum of the coordinates.
We omit the remaining details and just give the summary
of this step (for simplicity we ignore the constant $\lambda_k$):

\begin{quotation}
\noindent
Step ${\rm B}_k$. Given $x:=(x_1, \ldots, x_k) \in \sP_k$. 
Let $s := \sum_{i=1}^{k} x_i$, $r := \lfloor \sqrt{k} s \rfloor$. \\
Compute $x' := x^{\phi_2^{-r}}$. \\
Pass $x''' := x' N_k$ to the next stage, and output 
$y_k := (\mbox{fractional~part~of~} \sqrt{k} s)/k$.
\end{quotation}

Since the number of computations needed at each step is linear, we conclude that:

\begin{theorem}\label{Th3}
Given $x \in \sO_n$, $\Theta(x) \in \Pi$ can
be computed in $O(n^2)$ steps.
\end{theorem}

\noindent{\bf Remarks.}
The inverse map $\Theta ^{-1}$ is just as easy to compute,
since each of the individual steps is easily reversed.
Two details are worth mentioning.
When inverting step ${\rm B}_k$, given $x'''$ and $y_k$, we obtain $x'$ by
multiplying $x'''$ by $N_k ^{\tr}$
and adding $y_k/\sqrt{k}$ to each component.
For the computation of $r$,
it can be shown (we omit the proof) that
for inverting step ${\rm A}$, to go from $x'$ to $x$,
$r$ should be taken to be the number of strictly negative
components in $x'$. For step ${\rm B}_k$, $r$ is
the number of indices $i$ such that
$$
%\frac{x_i'}{\sqrt{k}} ~<~ b \, \sum_{i=1}^{k} x_i' \,.
x_i'  ~<~ b \, \sqrt{k} \, \sum_{j=1}^{k} x_j' \,.
$$

%%%% SECTION 5
\section{An alternative dissection of $\sO_4$}\label{Sec5}

In this section we give a six-piece dissection of $\sT := \sO_4$
into a prism $\sqrt{ \frac{3}{4} } \sP_3 \times I_{\frac{1}{2}}$.
Although it requires two more pieces than the dissection of
Theorem \ref{Th2}, it still only
uses three cuts.
It also has an appealing symmetry.

\begin{figure}[htb]
\begin{center}
\includegraphics[width=7.5cm]{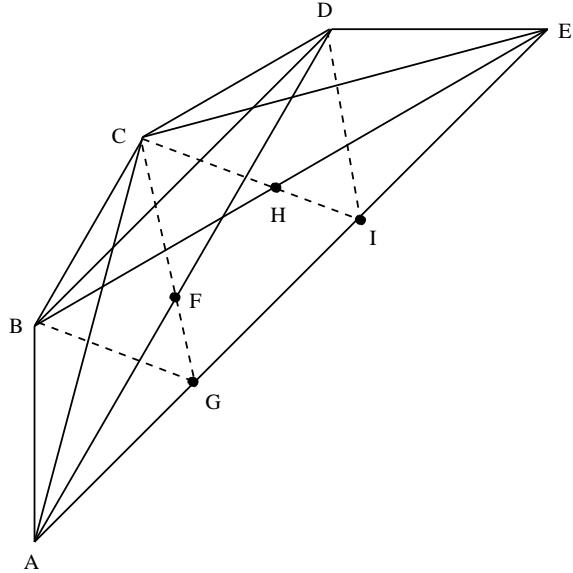}
\end{center}
\caption{ 
$\sT := \sO_4$ is the convex hull
of $A,B,C,D,E$;
the first two cuts are made along the hyperplanes containing
$B,C,F,G$ and $C,D,H,I$, respectively. }
\label{Fig5a}
\end{figure}

We start by subtracting $\frac{1}{2}$ from the coordinates in \eqn{EqH1},
in order to move the origin to the center of $\sT$.
That is, we take $\sT$ to be the convex hull of the points
$A:=( -\frac{1}{2}, -\frac{1}{2}, -\frac{1}{2}, -\frac{1}{2})$, 
$B:=( \frac{1}{2}, -\frac{1}{2}, -\frac{1}{2}, -\frac{1}{2})$, 
$C:=( \frac{1}{2}, \frac{1}{2}, -\frac{1}{2}, -\frac{1}{2})$, 
$D:=( \frac{1}{2}, \frac{1}{2}, \frac{1}{2}, -\frac{1}{2})$, 
$E:=( \frac{1}{2}, \frac{1}{2}, \frac{1}{2}, \frac{1}{2})$
(see Fig. \ref{Fig5a}). We use $(w,x,y,z)$ for coordinates in $\RR^4$.
Note that $\sT$ is fixed by the symmetry $(w,x,y,z) \mapsto
(-z,-y,-x,-w)$. 

%% This is Fig 5a

We make two initial cuts, along the hyperplanes $w+y+z = -\frac{1}{2}$
and $w+x+z = \frac{1}{2}$.  The first cut intersects the
edges of $\sT$ at the points $B$, $C$, $F:=(0,0,0,-\frac{1}{2})$
and $G:=( -\frac{1}{6}, -\frac{1}{6}, -\frac{1}{6}, -\frac{1}{6})$;
the second at the points $C$, $D$, $H:=(\frac{1}{2},0,0,0)$
and $I:=( \frac{1}{6}, \frac{1}{6}, \frac{1}{6}, \frac{1}{6})$.
The three pieces resulting from these cuts will
be denoted by $\sT_1$ (containing $A$),
$\sT_2$ (the central piece),
and $\sT_3$ (containing $E$).

We apply the transformation
$\alpha := (w,x,y,z) \mapsto (-y,-x,-w,-1-z)$ to $\sT_1$
and
$\beta := (w,x,y,z) \mapsto (1-w,-z,-y,-x)$ to $\sT_3$.
$\alpha$ fixes the triangle $BCF$, although not
pointwise, and similarly $\beta$ fixes the triangle $CDH$,
and so these transformations may be regarded as hinged,
in a loose sense of that word\footnote{Since
a hinged rod in the plane has a fixed point,
and a hinged door in three dimensions
has a fixed one-dimensional subspace,
a hinged transformation in
four dimensions should, strictly speaking,
have a two-dimensional region that is {\em pointwise} fixed.}.
This is what led us to
this dissection---we were attempting
to generalize Sch\"obi's hinged three-dimensional
dissection.

\begin{figure}[htb]
\begin{center}
\includegraphics[width=7.5cm]{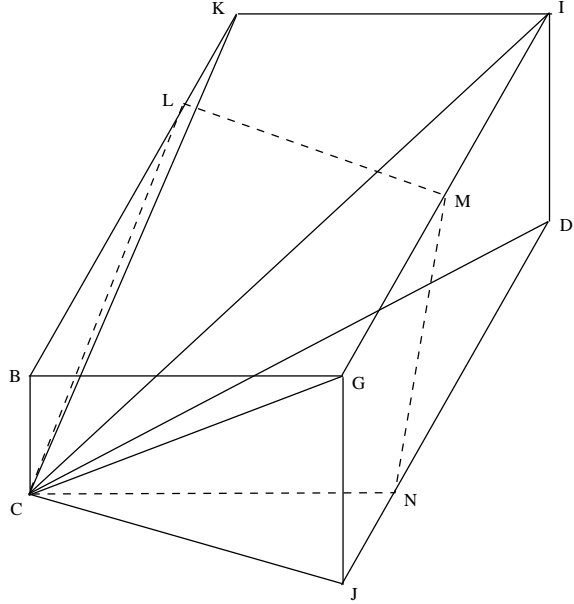}
\end{center}
\caption{ 
After the first two motions, we
have a polytope $\sT_4$ with seven vertices $B,G,J,C,K,I,D$.
The third cut is along the hyperplane containing $C,L,M,N$. }
\label{Fig5b}
\end{figure}

After applying $\alpha$ and $\beta$,
the resulting polytope $\sT_4 := \sT_1^{\alpha}+\sT_2+\sT_3^{\beta}$
is a convex body with seven vertices and six faces. (This
and other assertions in
this section were verified with the help of the 
programs Qhull~\cite{qhull}
and MATLAB~\cite{MATLAB}.) 
The seven vertices are $B, C, D, G, I$, 
$J := ( \frac{1}{6}, \frac{1}{6}, \frac{1}{6}, -\frac{5}{6})$
and
$K := ( \frac{5}{6}, -\frac{1}{6}, -\frac{1}{6}, -\frac{1}{6})$,
which are shown schematically in Fig. \ref{Fig5b}.
This figure is realistic in so far as it suggests
that the edges $BK$, $GI$ and $JD$ are equal
and parallel (in fact, $K-B= I-G = D-J =
( \frac{1}{3}, \frac{1}{3}, \frac{1}{3}, \frac{1}{3})$). 

%% This is Fig 5b

We now make one further cut, along the
hyperplane $w+x+y+z=0$, which 
separates $\sT_4$ into two pieces $\sT_5$ (containing $B$)
and $\sT_6$ (containing $K$).
This hyperplane meets the edge $BK$ at the point
$L:=( \frac{3}{4}, -\frac{1}{4}, -\frac{1}{4}, -\frac{1}{4})$,
$GI$ at the point
$M:=( 0, 0, 0, 0)$,
and $JD$ at the point
$N:=( \frac{1}{4}, \frac{1}{4}, \frac{1}{4}, -\frac{3}{4})$.
The point $L$ is three-quarters of the way along $BK$,
$M$ bisects $GI$, and
$N$ is one-quarter of the way along $JD$,

The final motion is to apply
$\gamma := (w,x,y,z) \mapsto (x,y,z,w-1)$
to $\sT_6$, and to form $\sT_7 := \sT_5 + \sT_6^{\gamma}$.
The convex hull of $\sT_7$ involves three new points,
the images of $L$, $M$ and $N$ under $\gamma$,
namely 
$P := ( -\frac{1}{4}, -\frac{1}{4}, -\frac{1}{4}, -\frac{1}{4})$,
$Q:=( 0, 0, 0, -1)$ and
$R := ( \frac{1}{4}, \frac{1}{4}, -\frac{3}{4}, -\frac{3}{4})$,
respectively.
Then $\sT_7$ is the convex hull of the eight points
$C, L, M, N$ and $R, B, P, Q$,
and it may be verified that the first four and the last
four of these points define copies of $\sqrt{ \frac{3}{4} } \sP_3$,
and that $\sT_7$ is
indeed congruent to $\sqrt{ \frac{3}{4} } \sP_3 \times I_{\frac{1}{2}}$,
as claimed.

We end with a question: can this construction be generalized to higher
dimensions?

\vspace*{+.1in}

\noindent{\bf Acknowledgments.}
We thank G.~N.~Frederickson for comments
on the manuscript and for drawing our attention to reference \cite{Eves1966}.

%%%BBBBBBB

\end{document}